\documentclass{amsart}
\usepackage{amsmath,amsthm,amssymb,verbatim}

\newtheorem{theorem}{Theorem}
\newtheorem*{mainthm}{Theorem \ref{thm:main}} 
\newtheorem*{thmtoric}{Theorem \ref{thm:toric}}
\newtheorem*{thmRX}{Theorem \ref{thm:RX}}
\newtheorem{lemma}[theorem]{Lemma}
\newtheorem{conj}[theorem]{Conjecture}
\newtheorem{cor}[theorem]{Corollary}
\newtheorem*{corss}{Corollary \ref{cor:ss}}

\begin{document}

\author{Ovidiu Munteanu and G\'abor Sz\'ekelyhidi}
\title{On convergence of the K\"ahler-Ricci flow}
\date{}

\begin{abstract}
We study the convergence of the K\"ahler-Ricci flow on a Fano manifold
under some stability conditions. More precisely we assume that the first
eingenvalue of the $\bar{\partial}$-operator
acting on vector fields is uniformly
bounded along the flow, and in addition the Mabuchi energy decays at
most logarithmically. 
We then give different situations in which the condition on the Mabuchi
energy holds. 
\end{abstract}

\maketitle

\section{Introduction}
Let $X$ be a compact K\"ahler manifold of dimension $n$ with $c_1\left(
X\right) >0$ and let us consider the Ricci flow introduced by
Hamilton~\cite{Ham82} defined by 
\[ \begin{aligned}
\frac \partial {\partial t}g_{i\bar{j}} & =-\left( R_{i\bar{j}}-g_{i\bar{j} 
}\right) \\
g_{i\bar{j}}\left( 0\right) &\in c_1(X).
\end{aligned}\]
Cao~\cite{Cao85} showed that this flow exists for all time. When
$c_1(X)=0$ or $c_1(X) < 0$, then Cao also showed, using Yau's
estimates~\cite{Yau78}, that the analogous normalized flow converges to a
K\"ahler-Einstein metric on $X$. When $c_1(X)>0$, ie. the manifold is
Fano,  then it is still an
open problem to determine whether $X$ admits a K\"ahler-Einstein metric.
The central problem in the field is the following conjecture.
\begin{conj}[Yau-Tian-Donaldson Conjecture] A Fano manifold $X$ admits a
	K\"ahler-Einstein metric if and only if it is K-polystable.
\end{conj}
For more details see \cite{Tian97}, \cite{Don02} and also \cite{PS08} for
a survey and many more references. If we try to use the K\"ahler-Ricci
flow to find a K\"ahler-Einstein metric then in light of this
conjecture the key problem is to relate K-polystability of $X$ to
convergence of the flow. This still seems out of reach at present, but
many partial results have been obtained. One possibility is to assume
that the Riemann curvature tensor remains uniformly bounded along the
flow, for results in this direction see \cite{PS06_1}, \cite{PSSW07_2},
\cite{GSz08}, \cite{To08}. 

In this paper we study a different kind of
assumption, which was introduced in \cite{PSSW07}. The main assumption
is that along the flow $g(t)$ there is a uniform lower bound $\lambda_t
\geqslant \lambda > 0$ on the lowest positive eigenvalue $\lambda_t$ of
the $\bar{\partial}$-operator 
acting on vector fields. The main result in \cite{PSSW07}
is that if in addition we assume that the Mabuchi energy is bounded from
below, then the flow converges to a K\"ahler-Einstein metric. Our main
result is a weakening of the hypothesis on the Mabuchi energy.

\begin{mainthm} 
	Suppose that along the K\"ahler-Ricci flow we have a
	uniform lower bound $\lambda_t\geqslant\lambda > 0$ on the first
	eingenvalue, and in addition the Mabuchi energy satisfies
	\begin{equation}\label{eq:Mlower}
		\mathcal{M}(g(t)) > -C\log(1+t) - D,
	\end{equation}
	for some constants $C,D>0$. Then the metrics $g(t)$ converge
	exponentially fast in $C^\infty$ to a K\"ahler-Einstein metric. 
\end{mainthm}

\noindent 
In fact the proof of the theorem shows that instead of (\ref{eq:Mlower})
it is enough to assume 
\begin{equation}\label{eq:infY}
	\inf_{t>0} Y(t)=0,
\end{equation}
\noindent where
\[
	Y(t) = -\frac{d}{dt}\mathcal{M}(g(t)).
\]
Since the Mabuchi energy is monotonically decreasing under the
K\"ahler-Ricci flow, the assumption of the theorem clearly implies this
weaker statement. We will prove this theorem in Section \ref{sec:main}.
The main novelty is the estimate in Lemma \ref{lem:estimate}.  

The advantange of replacing the lower bound on $\mathcal{M}$ by these
weaker statements is that there are some other natural conditions under
which they can be shown to hold. We give two such conditions. 

In Section
\ref{sec:toric} we show directly that on a K-semistable toric
variety the condition (\ref{eq:Mlower})
of Theorem \ref{thm:main} holds. 
\begin{thmtoric}
	Suppose that the Fano toric variety $X$ is K-semistable. 
	If $g(t)$ are torus invariant metrics satisfying the
	K\"ahler-Ricci flow then 
	\[ \mathcal{M}(g(t)) > -C\log(1+t) - D\]
	for some constants $C,D>0$. 
\end{thmtoric}	
\noindent While it is
already known that K-semistable toric varieties admit K\"ahler-Einstein
metrics (see \cite{WZ04}) and hence the Mabuchi energy is bounded below, 
this direct argument could be of independent interest. 

In Section \ref{sec:RM} we show that if
\[ R(X):=\sup\{t\,|\, \text{there exists } \omega\in c_1(X) \text{ such that }
\mathrm{Ric}(\omega) > t\omega\} =1,\]
then the condition (\ref{eq:infY}) holds.
\begin{thmRX}
	Suppose that $X$ is a Fano manifold which satisfies $R(X)=1$.
	Then along the K\"ahler-Ricci flow $g(t)$ on $X$, we have
	\[ \inf_{t > 0} Y(g(t)) = 0.\]
\end{thmRX}
\noindent As a corollary we note 
\begin{corss}
	If $X$ is a Fano manifold with $R(X)=1$ then $X$ is
	K-semistable.
\end{corss}
\noindent This is a strengthening of a result in \cite{GSz09}, where
the second author conjectured that the condition
$R(X)=1$ is equivalent to K-semistability. 

Finally we note that if $X$ is a degree $n+1$ hypersurface in
$\mathbf{P}^{n+1}$ and $n+1\geqslant 3$, then the $\alpha$-invariant of
$X$, defined by Tian~\cite{Tian87}, satisfies (see \cite{CP02})
\[ \alpha(X) \geqslant \frac{n}{n+1}. \]
This implies that $R(X)=1$, since given any K\"ahler metric $\omega_0$ 
and $0 < t < 1$ we can find $\omega$ that satisfies
\[ \mathrm{Ric}(\omega) = t\omega + (1-t)\omega_0,\]
so in particular $\mathrm{Ric}(\omega) > t\omega$ (see \cite{Tian87}).
By the above results
$X$ is K-semistable and 
also the condition (\ref{eq:infY}) holds. Note that it is not known
whether the Mabuchi energy is bounded from below for such hypersurfaces.

\subsection*{Acknowledgements}
We would like to thank D. H. Phong for his encouragement and interest in
this work. 

\section{Background}
In this section we recall some basic notation and results that we use. 
We consider the normalized K\"ahler-Ricci flow
\[ \begin{aligned}
	\frac{\partial}{\partial t} g_{i\bar{j}} &= - (R_{i\bar{j}} -
	g_{i\bar{j}}) = \partial_i\partial_{\bar{j}} u,\\
	g_{i\bar{j}}(0) &\in c_1(X). 
\end{aligned}\]
Here $R_{i\bar{j}}$ is the Ricci curvature and
$u(t)$ is the Ricci potential of the metric $g(t)$, which we
normalize so that 
\[ \frac{1}{V}\int_X e^{-u}\omega^n = 1.\]
Here 
\[ V = \int_X \omega^n.\]
A fundamental result of Perelman \cite{Per02} (see \cite{ST03} for a
detailed exposition) is that there exists a
constant $C_0$ depending only on $g(0)$ such that 
\begin{equation}\label{eq:perelman}
	\Vert u\Vert_{C^0} + \Vert \nabla u\Vert_{C^0} + \Vert \Delta
u\Vert_{C^0} \leqslant C_0
\end{equation}
along the flow. We also need to recall the following Lemma from
\cite{PSSW07}. 
\begin{lemma}\label{lem:smoothing} 
	We have the following two results along the K\"ahler-Ricci flow.

	\begin{itemize}
		\item There exists $\delta, K>0$ depending on the
			dimension $n$ with the following property. For
			any $\epsilon$ with $0 < \epsilon \leqslant
			\delta$ and any $t_0\geqslant 0$, if
			\[\Vert u(t_0)\Vert_{C^0} \leqslant\epsilon,\]
			then
			\[\Vert \nabla u(t_0+2)\Vert_{C^0} + \Vert
			R(t_0+2)-n\Vert_{C^0} \leqslant K\epsilon.
			\]
		\item There exists $C>0$ depending on $g(0)$ such that
			\[ \Vert u\Vert_{C^0}^{n+1}\leqslant C\Vert
			\nabla u\Vert_{L^2}\Vert \nabla
			u\Vert_{C^0}^n,\]
	\end{itemize}
\end{lemma}
\begin{proof}
	The first part is Lemma 1 from \cite{PSSW07}. The second part also
	follows directly from Lemma 3 in \cite{PSSW07}. There it is shown
	that 
	\[ \Vert u-b\Vert_{C^0}^{n+1}\leqslant C\Vert
			\nabla u\Vert_{L^2}\Vert \nabla
			u\Vert_{C^0}^n,\]
	where $b$ is the average
	\[ b = \frac{1}{V}\int_X ue^{-u}\,\omega^n,\]
	and also $|b|\leqslant \Vert u-b\Vert_{C^0}$. But then
	\[ \Vert u\Vert_{C^0}\leqslant \Vert u-b\Vert_{C^0} +
	|b|\leqslant 2\Vert u-b\Vert_{C^0},\]
	so our statement follows. 
\end{proof}

We define the Mabuchi functional $\mathcal{M}$ so that
$\mathcal{M}(g(0))=0$ and 
\begin{equation}\label{eq:Mab}
	\frac{d}{dt}\mathcal{M}(g(t)) = -\int_X u(R - n)\,\omega^n = -\int_X
|\nabla u|^2\,\omega^n,
\end{equation}
where $R$ is the scalar curvature of the metric $g(t)$. In particular
$\mathcal{M}$ is monotonically decreasing under the flow. We define
\[ Y(t) = \int_X |\nabla u|^2\,\omega^n.\]
The key to proving convergence of the K\"ahler-Ricci flow is showing
exponential decay of $Y(t)$ (see \cite{PSSW07}). For this the basic
inequality proved in \cite{PS06_1} is
\begin{equation}\label{1}
	\begin{aligned} 
\frac d{dt}Y\left( t\right)  \leqslant& -2\lambda _tY\left( t\right) -2\lambda
_tFut\left( \pi _t\left( \nabla ^ju\right) \right) -\int_X\left| \nabla
u\right| ^2\left( R-n\right) \omega ^n   \\
&-\int_X\nabla ^ju\nabla ^{\bar{k}}u\left( R_{j\bar{k}}-g_{j\bar{k}}\right)
\omega ^n,  
	\end{aligned}
\end{equation}
where $Fut\left( \pi _t\left( \nabla ^ju\right) \right) $ is the Futaki
invariant, applied to the orthogonal projection $\pi _t\left( \nabla
^ju\right) $ of the vector field $\nabla ^ju$ on the space of holomorphic
vector fields.

\section{The main argument}\label{sec:main}
Our goal is to prove the following result.

\begin{theorem}\label{thm:main}
Assume that the lowest positive eigenvalue of the Laplacian $-g^{j\bar{k}%
}\nabla _j\nabla _{\bar{k}}$ acting on $T^{1,0}$ vector fields has a
positive lower bound and that the Mabuchi energy along the
K\"ahler-Ricci flow satisfies
\[ \mathcal{M}(g(t)) > -C\log(1+t) - D,\]
for some constants $C,D>0$. 
Then the metrics $g_{i\bar{j}}$
converge exponentially fast in $C^\infty $ to a K\"ahler-Einstein metric.
\end{theorem}

Before proving the theorem we need two Lemmas. 
\begin{lemma}\label{lem:Futaki}
	If $\inf_{t\geqslant 0} Y(t)=0$, then the Futaki invariant of
	$X$ vanishes.
\end{lemma}
\begin{proof}
	Using the second part of Lemma \ref{lem:smoothing} we have
\[
\Vert u\Vert _{C^0}^{n+1}\left( t\right) \leqslant C\Vert \nabla u\Vert
_{L^2}\left( t\right) \Vert \nabla u\Vert _{C^0}^n\left( t\right)
\leqslant C_1 Y(t)^{1/2}, 
\]
where we have also used Perelman's estimate (\ref{eq:perelman}). 
It follows that
$\inf_{t\geqslant 0}\Vert u\Vert _{C^0}= 0.$
Then the first part of Lemma \ref{lem:smoothing} implies
that 
\[
\inf_{t\geqslant 0}\Vert R-n\Vert _{C^0}=0.
\]
Thus the manifold has to be K-semistable by Donaldson's lower bound
\cite{Don05}. In particular the Futaki invariant of $X$ vanishes.
\end{proof}

\begin{lemma}\label{lem:estimate}
On any K\"ahler manifold $X$ such that 
\[
R_{i\bar{j}}-g_{i\bar{j}}=-u_{i\bar{j}}
\]
we have the following estimate 
\[
\left| \int_Xu_{j\bar{k}}\nabla ^ju\nabla ^{\bar{k}}u\right| \leqslant 5\left(
\Vert \nabla u\Vert _{C^0}^2+\Vert \Delta u\Vert _{C^0}\right) \int_X\left|
\nabla u\right| ^2
\]
\end{lemma}
\begin{proof}
We have 
\begin{equation} \label{eq:1}
	\begin{aligned}
\left| \int_Xu_{j\bar{k}}\nabla ^ju\nabla ^{\bar{k}}u\right| \leqslant
&\int_X\left| u_{j\bar{k}}\right| \left| \nabla u\right| ^2 \\
\leqslant &\left( \int_X\left| u_{j\bar{k}}\right| ^2\left| \nabla u\right|
^2\right) ^{\frac 12}\left( \int_X\left| \nabla u\right| ^2\right) ^{\frac
12}.
\end{aligned}
\end{equation}
We now denote 
\[
I=\int_X\left| u_{j\bar{k}}\right| ^2\left| \nabla u\right| ^2. 
\]
Integration by parts yields: 
\begin{equation}\label{3}
	\begin{aligned}
I &=\int_Xu_{j\bar{k}}u_{k\bar{j}}\left| \nabla u\right| ^2  \\
&=-\int_Xu_{j\bar{k}\bar{j}}u_k\left| \nabla u\right| ^2-\int_Xu_{j\bar{k}
}u_k\left( \left| \nabla u\right| ^2\right) _{\bar{j}} \\
&\leqslant \left| \int_Xu_{j\bar{k}\bar{j}}u_k\left| \nabla
u\right|^2\right| + \left|  \int_Xu_{j\bar{k}
}u_k\left( \left| \nabla u\right| ^2\right) _{\bar{j}}  \right|. 
	\end{aligned}
\end{equation}
The first term above is, using the Ricci identities, 
\[\begin{aligned}
\left|\int_Xu_{j\bar{k}\bar{j}}u_k\left| \nabla u\right| ^2\right|
&=\left|\int_X\left(
\Delta u\right) _{\bar{k}}u_k\left| \nabla u\right| ^2\right| \\
&\leqslant\int_X\left( \Delta u\right) ^2\left| \nabla u\right|^2+
\left|\int_X 
(\Delta u) u_k\left( \left| \nabla u\right| ^2\right) _{\bar{k}}\right| \\
&\leqslant \Vert \Delta u\Vert _{C^0}^2\int_X\left| \nabla u\right|
^2+\int_X\left| \Delta u\right| \left| \nabla u\right| \left| \nabla \left(
\left| \nabla u\right| ^2\right) \right| \\
&\leqslant \Vert \Delta u\Vert _{C^0}^2\int_X\left| \nabla u\right| ^2+\frac
12\int_X\left( \Delta u\right) ^2\left| \nabla u\right| ^2+\frac
12\int_X\left| \nabla \left( \left| \nabla u\right| ^2\right) \right| ^2 \\
&\leqslant \frac 32\Vert \Delta u\Vert _{C^0}^2\int_X\left| \nabla u\right|
^2+\frac 12\int_X\left| \nabla \left( \left| \nabla u\right| ^2\right)
\right| ^2.
\end{aligned}\]

\noindent The second term can be estimated as follows: 
\[\begin{aligned}
\left|\int_Xu_{j\bar{k}}u_k\left( \left| \nabla u\right| ^2\right)
_{\bar{j}}\right|
&\leqslant \int_X\left| u_{j\bar{k}}\right| \left| \nabla u\right| \left| \nabla
\left( \left| \nabla u\right| ^2\right) \right| \\
&\leqslant \frac 12\int_X\left| u_{j\bar{k}}\right| ^2\left| \nabla u\right|
^2+\frac 12\int_X\left| \nabla \left( \left| \nabla u\right| ^2\right)
\right| ^2 \\
&=\frac 12I+\frac 12\int_X\left| \nabla \left( \left| \nabla u\right|
^2\right) \right| ^2.
\end{aligned}\]

\noindent Using these estimates in (\ref{3}) it follows that 
\begin{equation}
I\leqslant 3\Vert \Delta u\Vert _{C^0}^2\int_X\left| \nabla u\right|
^2+2\int_X\left| \nabla \left( \left| \nabla u\right| ^2\right) \right| ^2.
\label{4}
\end{equation}

\noindent We now denote 
\[\begin{aligned}
J &=\int_X\left| \nabla \left( \left| \nabla u\right| ^2\right) \right| ^2
\\
&=-\int_X\left| \nabla u\right| ^2\Delta \left| \nabla u\right| ^2.
\end{aligned}\]

\noindent According to the Bochner formula, 
\[
\Delta \left| \nabla u\right| ^2=2\langle\nabla u,\nabla \left( \Delta u\right)
\rangle+R_{i\bar{j}}u_{\bar{\imath}}u_j+\left| u_{i\bar{j}}\right| ^2+\left|
u_{ij}\right| ^2. 
\]

\noindent Notice that 
\[\begin{aligned}
R_{i\bar{j}}u_{\bar{\imath}}u_j &=\left( g_{i\bar{j}}-u_{i\bar{j}}\right)
u_{\bar{\imath}}u_j \\
&\geqslant -u_{i\bar{j}}u_{\bar{\imath}}u_j \\
&\geqslant -\left| u_{i\bar{j}}\right| ^2-\frac 14\left| \nabla u\right| ^4,
\end{aligned}\]

\noindent where the last inequality follows from 
\[\begin{aligned}
\left| u_{i\bar{j}}u_{\bar{\imath}}u_j\right| &\leqslant \left| u_{i\bar{j}%
}\right| \left| \nabla u\right| ^2 \\
&\leqslant \left| u_{i\bar{j}}\right| ^2+\frac 14\left| \nabla u\right| ^4.
\end{aligned}\]

\noindent This proves that 
\[
\Delta \left| \nabla u\right| ^2\geqslant 2\langle\nabla u,\nabla \left( \Delta
u\right) \rangle-\frac 14\left| \nabla u\right| ^4. 
\]

\noindent We use this to estimate $J$ from above: 
\begin{equation}\label{5}
	\begin{aligned}
J &=-\int_X\left| \nabla u\right| ^2\Delta \left| \nabla u\right| ^2 
 \\
&\leqslant -2\int_X\langle\nabla u,\nabla \left( \Delta u\right)
\rangle\left| \nabla
u\right| ^2+\frac 14\int_X\left| \nabla u\right| ^4\left| \nabla u\right| ^2.
	\end{aligned}
\end{equation}

\noindent Let us study the first term in (\ref{5}). We have 
\[\begin{aligned}
-2\int_X\langle\nabla u,\nabla \left( \Delta u\right) &\rangle
\left| \nabla u\right|
^2= 2\int_X\left( \Delta u\right) ^2\left| \nabla u\right| ^2 
+2\int_X\left( \Delta u\right) \left\langle\nabla u,\nabla \left( \left| \nabla
u\right| ^2\right) \right\rangle \\
&\leqslant 2\Vert \Delta u\Vert _{C^0}^2\int_X\left| \nabla u\right|
^2+2\int_X\left| \Delta u\right| \left| \nabla u\right| \left| \nabla \left(
\left| \nabla u\right| ^2\right) \right| \\
&\leqslant 2\Vert \Delta u\Vert _{C^0}^2\int_X\left| \nabla u\right|
^2+2\int_X\left| \Delta u\right| ^2\left| \nabla u\right| ^2+\frac
12\int_X\left| \nabla \left( \left| \nabla u\right| ^2\right) \right| ^2 \\
&\leqslant 4\Vert \Delta u\Vert _{C^0}^2\int_X\left| \nabla u\right| ^2+\frac
12J.
\end{aligned}\]

\noindent Plugging this into (\ref{5}) we get 
\[
J\leqslant 8\Vert \Delta u\Vert _{C^0}^2\int_X\left| \nabla u\right| ^2+\frac
12\Vert \nabla u\Vert _{C^0}^4\int_X\left| \nabla u\right| ^2. 
\]

\noindent We plug this into (\ref{4}) and obtain 
\[
I\leqslant \left( 19\Vert \Delta u\Vert _{C^0}^2+\Vert \nabla u\Vert
_{C^0}^4\right) \int_X\left| \nabla u\right| ^2. 
\]
Using this in (\ref{eq:1}) we obtain the result.
\end{proof}

We can now prove the theorem.
\begin{proof}[Proof of Theorem~\ref{thm:main}]
Note first of all that our hypothesis implies that 
\begin{equation}\label{eq:inf}
	\inf_{t\geqslant 0} Y(t) = 0.
\end{equation}
For if $Y(t) > \epsilon>0$ for all $t$ then by Equation (\ref{eq:Mab})
\[ \mathcal{M}(g(t)) < -\epsilon t, \]
which contradicts our assumption.

Our goal is to prove that $Y\left( t\right) $ has exponential decay, since
then the exponential convergence of the K\"ahler-Ricci flow follows like in 
\cite{PSSW07}.
In the Inequality (\ref{1}) we use our
hypothesis that $\lambda_t\geqslant \lambda > 0$ and that the Futaki
invariant vanishes by Lemma \ref{lem:Futaki}. We obtain
\begin{equation}
\frac d{dt}Y\left( t\right) \leqslant -2\lambda Y\left( t\right) -\int_X\left|
\nabla u\right| ^2\left( R-n\right) \omega ^n-\int_X\nabla ^ju\nabla ^{\bar{k%
}}u\left( R_{j\bar{k}}-g_{j\bar{k}}\right) \omega ^n.  \label{2}
\end{equation}
Now using Lemma \ref{lem:estimate} and the fact that $R-n=-\Delta u$ we get
\begin{equation}
\frac d{dt}Y\left( t\right) \leqslant -2\lambda Y\left( t\right) +6\left( \Vert
\Delta u\Vert _{C^0}\left( t\right) +\Vert \nabla u\Vert _{C^0}^2\left(
t\right) \right) Y\left( t\right) .  \label{6}
\end{equation}

\noindent We remark that formula (\ref{6}) can be used as a substitute
for the
differential-difference inequality (5.5) in \cite{PSSW07}. 

We are now ready to finish the proof of the Theorem.
Fix $\varepsilon _0>0$ small to be determined later.
There must exist a point $t_0>0$ such that 
\[
Y\left( t_0\right) \leqslant \varepsilon _0
\]
because of (\ref{eq:inf}). It follows from
(\ref{6}) and Perelman's estimates on $\nabla u$ and $\Delta u$ that
\[
\frac d{dt}Y\left( t\right) \leqslant CY\left( t\right)
\]
for some constant $C>0$. 
Hence $Y$ has at most exponential growth and it follows that 
\[
Y\left( t_0+2\right) \leqslant Y\left( t_0\right) e^{2C}\leqslant \varepsilon _0e^{2C}.
\]
Consequently, if we set 
\[
\varepsilon _1=2\varepsilon _0e^{2C}, 
\]
then for $t\in\left[ t_0,t_0+2\right] $ we have 
\begin{equation}\label{Y}
Y\left( t\right) \leqslant \frac 12\varepsilon _1. 
\end{equation}

Assume that there exists a time $t_0<t<\infty $ for which $Y\left( t\right)
>\varepsilon _1.$ Then let
\[ t_1 = \inf\{t\,|\, t > t_0,\mbox{ and }Y(t)=\varepsilon_1\}\]
be the first time after $t_0$ such that $Y\left( t_1\right)
=\varepsilon_1$. By (\ref{Y}) we have $t_1 > t_0+2$ and so by the
definition of $t_1$ we have 
$Y\left( t_1-2\right) <\varepsilon_1$. By the second part of Lemma
\ref{lem:smoothing}  
\[
\Vert u\Vert _{C^0}\left( t_1-2\right) \leqslant C\varepsilon _1^{1/2\left(
n+1\right) }.
\]
Moreover if $\varepsilon_1$ is sufficiently small,
the first part of Lemma \ref{lem:smoothing} gives that 
\[
\Vert \Delta u\Vert _{C^0}\left( t_1\right) +\Vert \nabla u\Vert
_{C^0}^2\left( t_1\right) \leqslant C'\varepsilon _1^{1/2\left( n+1\right) }.
\]
Denote 
\[
\varepsilon =6C'\varepsilon _1^{1/2\left( n+1\right) },
\]
then (\ref{6}) implies that 
\[
\left( \frac d{dt}Y\right) \left( t_1\right) \leqslant -2\lambda Y\left(
t_1\right) +\varepsilon Y\left( t_1\right) . 
\]

Choose $\varepsilon_0$ small enough so that
$\varepsilon<\lambda$. 
Then we have 
\[
\left( \frac d{dt}Y\right) \left( t_1\right) \leqslant -\lambda Y\left(
t_1\right) =-\lambda \varepsilon _1<0, 
\]
which shows that $Y$ is decreasing in a neighborhood of $t_1$, and
contradicts the choice of $t_1$.
The contradiction came from our assumption that there exists finite $t$ such
that $Y\left( t\right) >\varepsilon_1$. Therefore, for all $t>t_0$ we
must have 
\[
Y\left( t\right) \leqslant\varepsilon_1. 
\]

Then as above for any $t\geqslant t_0+2$ we have 
\[
\left( \frac d{dt}Y\right) \left( t\right) \leqslant -\lambda Y\left( t\right) .
\]
This shows that $Y\left( t\right) $ is exponentially decreasing and
therefore by the argument in \cite{PSSW07} we get the exponential convergence
of the K\"ahler-Ricci flow to a K\"ahler-Einstein metric.
\end{proof}

\section{The case of toric varieties}\label{sec:toric}
In this section we prove the following.
\begin{theorem}\label{thm:toric}
	Suppose that the Fano toric variety $X$ is K-semistable. 
	If $g(t)$ are torus invariant metrics satisfying the
	K\"ahler-Ricci flow then 
	\[ \mathcal{M}(g(t)) > -C\log(1+t) - D\]
	for some constants $C,D>0$. 
\end{theorem}
In fact if $X$ is K-semistable, then its Futaki invariant must vanish,
so by Wang-Zhu \cite{WZ04}
$X$ admits a K\"ahler-Einstein metric. Then Bando-Mabuchi \cite{BM85}
implies that
the Mabuchi energy is bounded below. So a stronger result 
follows easily from known results, but the interest lies in our
more direct proof which uses K-semistability explicitly. 
The proof follows the argument in \cite{GSz07_2} for the
Calabi flow where also more details can be found. 
\begin{proof}
	Suppose that the torus invariant metrics
	$g(t)$ satisfy the K\"ahler-Ricci flow. On the dense complex
	torus $(\mathbf{C}^*)^n\subset X$ we have 
	$g(t) =
	i\partial\bar{\partial}\phi(t)$ for some torus invariant
	functions $\phi(t)$. We can therefore think of them as functions
	on Euclidean space:
	\[ \phi(t) : \mathbf{R}^n \to\mathbf{R}. \]
	These K\"ahler potentials $\phi(t)$ satisfy
	\[ \frac{\partial}{\partial t} \phi(t) = u(t) =
	\log\det(\phi_{ij}) + \phi,\]
	where $u$ is the Ricci potential as before. For each $t$ the
	symplectic potential $f(t)$ is the Legendre transform of $\phi(t)$.
	Then $f(t)$ is a convex function on a polytope $P$, satisfying
	the Guillemin boundary conditions (for more details see
	\cite{Don02}). We have
	\[ \frac{\partial}{\partial t} f = -L(f) =
	\log\det(f_{ij}) + f - x\cdot\nabla f, \]
	where $x$ is the Euclidean coordinate on the polytope $P$ and the
	function $L(f)$ is just the Ricci potential expressed in the $x$
	coordinates.
	Let $g$ be a fixed symplectic potential, and define the
	functional
	\[ \mathcal{F}(f) = -\int_P \log\det(g^{ik}f_{kj})\,d\mu +
	\int_P g^{ij}f_{ij}\,d\mu.\]
	Then 
	\[ \frac{d}{dt}\mathcal{F}(f) = \int_P f^{ij}L(f)_{ij}\,d\mu -
	\int_P g^{ij}L(f)_{ij}\,d\mu = \int_P (f^{ij} - g^{ij})
	L(f)_{ij}\,d\mu.\]
	When integrating by parts the boundary terms vanish, so
	\[ \frac{d}{dt}\mathcal{F}(f) = \int_P\left[(f^{ij})_{,ij} -
	(g^{ij})_{,ij}\right]L(f)\,d\mu = \int_P(R(g)-R(f))L(f)\,d\mu,\]
	where $R(f),R(g)$ are the scalar curvatures of the metrics
	determined by $f,g$. By Perelman's estimates $R(f)$ and the
	Ricci potential $L(f)$ (normalized by adding a constant)
	are uniformly bounded along the flow, so we obtain
	\[ \frac{d}{dt} \mathcal{F}(f(t)) < C_1,\]
	ie. 
	\begin{equation}\label{Fbound}
		\mathcal{F}(f(t)) < C_1t + C_2 
	\end{equation}
	for some constants $C_1, C_2$. 
	
	Applying the inequality $\log x <
	x/2$ to each eigenvalue, we obtain $\log\det(M)\leqslant
	\frac{1}{2}\mathrm{Tr}(M)$ 
	for any positive definite matrix $M$. Applying this to the
	defining formula of $\mathcal{F}$ we get
	\begin{equation}\label{Lbound}
		\mathcal{F}(f)\geqslant \frac{1}{2}\int_P
	g^{ij}f_{ij}\,d\mu.
	\end{equation}

	The AM-GM inequality implies that 
	\[ -\log\det(g^{ik}f_{kj}) \geqslant -n\log (g^{ij}f_{ij})\]
	so using the convexity of $-\log$ we get
	\begin{equation}\label{I}
		\begin{aligned}
		-\int_P\log\det(g^{ik}f_{kj})\,d\mu &\geqslant
		-n\int_P \log (g^{ij}f_{ij})\,d\mu  \\
		&\geqslant -C_3\log\int_P g^{ij}f_{ij}\,d\mu -C_4. \\
		&\geqslant -C_5\log(1+t) - C_6,
	\end{aligned}\end{equation}
	where in the last line we have used (\ref{Fbound}) and
	(\ref{Lbound}).

	In terms of symplectic potentials the Mabuchi energy is given by
	\[ \mathcal{M}(f) = -\int_P \log\det(f_{ij})\,d\mu +
	\int_{\partial P}f\,d\sigma - n\int_P f\,d\mu,\]
	moreover if the manifold is K-semistable then 
	\[ \int_{\partial P}f\,d\sigma - n\int_P f\,d\mu \geqslant 0\]
	for all convex functions $f$ (see Donaldson \cite{Don02}).
	Therefore we have
	\[ 
		\mathcal{M}(f)\geqslant -\int_P
		\log\det(g^{ik}f_{kj})\,d\mu -
		\int_P\log\det(g_{ij})\,d\mu \geqslant -C\log(1+t) -D
	\]
	using (\ref{I}). This completes the proof.
\end{proof}

\section{The case when $R(X)=1$}\label{sec:RM}
For a Fano manifold $X$ we define 
\[ R(X) = \sup\{t\, |\, \mbox{ there exists a metric } \omega\in c_1(X)
\mbox{ such that }\mathrm{Ric}(\omega) > t\omega\}.\]
We show the following.
\begin{theorem}\label{thm:RX}
	Suppose that $X$ is a Fano manifold which satisfies $R(X)=1$.
	Then along the K\"ahler-Ricci flow $g(t)$ on $X$, we have
	\[ \inf_{t > 0} Y(g(t)) = 0.\]
\end{theorem}
\begin{proof}
	We argue by contradiction. Suppose that 
	\[ \inf_{t>0} Y(g(t)) = \epsilon > 0.\]
	Since $\frac{d}{dt}\mathcal{M}(g(t)) = -Y(g(t))$, we then have
	\begin{equation}\label{eq:Mabuchi}
		\mathcal{M}(g(t)) < -\epsilon t + C
	\end{equation}
	for some constant $C$. Fix a base metric $\omega_0\in c_1(X)$,
	and define the $\mathcal{J}$ functional by
	$\mathcal{J}(\omega_0)=0$ and
	\[ \frac{d}{ds}\mathcal{J}(\omega_s) = \int_X \dot\phi_s
	(\Lambda_{\omega_s}\omega_0 - n)\,\omega_s^n,\]
	where $\omega_s=\omega_0+i\partial\bar{\partial}\phi_s$ is a
	path of metrics. Using the path $\omega_s = \omega_0 +
	si\partial\bar{\partial}\phi$ we get
	\[ \begin{aligned}
		\mathcal{J}(\omega_0 + i\partial\bar{\partial}\phi) &=
		\int_0^1\int_X \phi(\Lambda_{\omega_s}\omega_0 - n)\,
		\omega_s^n\,ds \\
		&= n\int_0^1\int_X \phi(\omega_0-\omega_s)\wedge
		\omega_s^{n-1}\,ds \\
		&= n\int_0^1\int_X \phi(-si\partial\bar{\partial}
		\phi)\wedge(s\omega_1 + (1-s)\omega_0)^{n-1}\,ds \\
		&=n\int_0^1\int_X \phi(-i\partial\bar{\partial}\phi)
		\wedge\sum_{k=0}^{n-1}\binom{n-1}{k}
		s^{k+1}(1-s)^{n-1-k}\omega_1^k\wedge\omega_0^{n-1-k}
		\,ds\\
		&= n\int_X \phi(-i\partial\bar{\partial}\phi) \wedge
		\sum_{k=0}^{n-1} \binom{n-1}{k}\frac{1}{n+1}
		\binom{n}{k+1}^{-1}\omega_1^k\wedge\omega_0^{n-1-k} \\
		&= \int_X \phi(-i\partial\bar{\partial}\phi)\wedge
		\sum_{k=0}^{n-1} \frac{n}{n+1}\cdot \frac{k+1}{n}
		\cdot\omega_1^k\wedge\omega_0^{n-1-k} \\
		&\leqslant\frac{n}{n+1}\int_X\phi(\omega_0-\omega_1)
		\wedge\sum_{k=0}^{n-1}\omega_1^k\wedge\omega_0^{n-1-k}\\
		&= \frac{n}{n+1}\int_X\phi(\omega_0^n-\omega_1^n).
	\end{aligned}\]
	This is the well-known inequality $I-J\leqslant\frac{n}{n+1}I$
	in the literature in terms of Aubin's $I,J$ functionals (see
	\cite{Aub84}). 
	
	The point is that along the K\"ahler-Ricci flow $g(t)$ we have
	$|\dot{\phi}|<C_1$ for some constant $C_1$ 
	by Perelman's estimates, so it follows that
	\[ |\phi(t)| < C_1t + C_2.\] 
	It follows that
	\[ \mathcal{J}(g(t)) < C_1t + C_2\]
	for some different constants $C_1,C_2$. But then using
	(\ref{eq:Mabuchi})
	\[ \mathcal{M}(g(t)) + \frac{\epsilon}{2C_1}\mathcal{J}(g(t)) <
	-\frac{\epsilon}{2}t + C',\]
	and in particular the functional $\mathcal{M} +
	\frac{\epsilon}{2C_1}\mathcal{J}$ is not bounded from below on
	$c_1(X)$. It follows then using the work of
	Chen-Tian~\cite{CT05_1} (see \cite{GSz09}) that 
	\[ R(X) \leqslant 1 - \frac{\epsilon}{2C_1},\]
	which is a contradiction.
\end{proof}

Finally we note
\begin{cor}\label{cor:ss}
	If $X$ is a Fano manifold with $R(X)=1$ then $X$ is
	K-semistable.
\end{cor}
\begin{proof}
	This follows from the previous theorem and the
	proof of Lemma \ref{lem:Futaki}. 
\end{proof}

\bibliographystyle{hacm}
\bibliography{../mybib}

\bigskip
\begin{flushright}
	{\sc Department of Mathematics \\
	Columbia University \\ New York, NY 10027}
\end{flushright}

\end{document}